\documentclass[reqno,12pt]{amsart}
\usepackage{graphicx}
\usepackage{epsfig}
\usepackage{amssymb,amsthm,amsmath,amsfonts,amstext,mathrsfs,fancybox}
\usepackage{enumerate}
\usepackage{asymptote}
\usepackage{latexsym}
\usepackage{epsfig}
\usepackage{url}
\usepackage{fancyhdr}
\usepackage{multirow}
\usepackage{pdfpages}
\usepackage{pdfpages}

\newtheorem*{thmm}{Theorem}

\newtheorem*{lemm}{Lemma}
\newtheorem{prop}{Proposition}
\newtheorem*{prob}{Problem}
\newtheorem*{defin}{Definition}

\makeatletter
\def\dual#1{\expandafter\dual@aux#1\@nil}
\def\dual@aux#1/#2\@nil{\begin{tabular}{@{}l@{}}#1\\#2\end{tabular}}
\makeatother

\input xy
\xyoption{all}

\title[Friendly paths]
{Friendly paths for finite subsets\\
of plane integer lattice. I}
\author[G. Alkauskas]{Giedrius Alkauskas}
\address{Vilnius University, Institute of Informatics, Naugarduko 24, LT-03225 Vilnius, Lithuania}
\email{giedrius.alkauskas@mif.vu.lt}
\date{\today}

\begin{document}
\begin{abstract}
For a given finite subset $\mathcal{P}$ of points of the lattice $\mathbb{Z}^{2}$, \emph{a friendly path} is a monotone (uphill or downhill) lattice path which splits points in half; points lying on the path itself are discarded. The purpose of this paper (and its sequel) is to fully describe all configurations of $n$ points in $\mathbb{Z}^2$ which do not admit a friendly path. We say that such an $n$-set is \emph{inseparable}. There are, up to the lattice symmetry, exactly $c(n)$ such sets. If only lattice shifts are counted, there are $\hat{c}(n)$ of them. Both sequences are new entries into OEIS (A369382 and, respectively, A367783). In particular, $n=27$ is the first odd numbers with $c(n)=1$. No example was known so far. This solves problem 11484(b)* posed in American Mathematical Monthly (February 2010). In this paper we also show that inseparable $n$-set exist for all even numbers $n\geq 12$ and almost all odd numbers.
\end{abstract}

\date{\today}

\maketitle
\section{Friendly paths}
\subsection{Formulation}The following problem appeared in American Mathematical Monthly (\cite{monthly1}, Problem 11484).

\begin{prob}An \emph{uphill} (NE, North-East) lattice path is the union of a (doubly infinite) sequence of directed line segments in $\mathbb{R}^{2}$, each connecting an integer pair $(a,b)$ to an adjacent pair, either $(a,b+1)$ or $(a+1,b)$. A \emph{downhill} (SE, South-East) lattice path is defined similarly, but with $b-1$ in place of $b+1$, and a \emph{monotone} lattice path is an uphill or downhill lattice path.\\

\indent Given a finite set $\mathcal{P}$ of points in $\mathbb{Z}^{2}$, a \emph{friendly} path is a monotone lattice path for which there are as many points in $\mathcal{P}$ on one side of the path as on the other (points that lie on the path do not count). 
\begin{itemize}
\item[(\textbf{a})\,\,\,]Prove that if $n=b^{2}+a^{2}+b+a$ for some positive integers $a,b$ such that $b\leq a\leq b+\sqrt{2b}$, then there exists a configuration of $n$ such points that there does not exist a friendly line.
\item[(\textbf{b})*]Is it true that for every odd-sized set of points there is a friendly path?
\end{itemize}
\end{prob}
Asterisk means that no solution was known to the presenter of the problem (myself). None was received by the editors in two years time \cite{monthly2}. \\

We remark that, as is clear from \cite{monthly2}, the precise bound for the problem (\textbf{a}) is slightly higher:
	\begin{eqnarray*}
		n=a^2+a+b^2+b,\,a,b\in\mathbb{N},\, b\leq a\leq b+\sqrt{2b}+\frac{1}{2}.
	\end{eqnarray*}
As is easy to prove, such representation of $n$, if it exists, is unique. The sequence of these numbers starts from
	\begin{eqnarray*}
		&&4, 8, 12, 18,24,26,32,40,42,50,60,62,72,76,84,86,98,102,112,114,128,132,144, 146, 162,\\
		&&166, 180,182,188,200,204,220,222,228,242,246,264,266,272,288,292,312,314,320,338,342.
	\end{eqnarray*}

\subsection{Results} In this paper we solve the problem and show that it further leads to more profound topics in geometric combinatorics and number theory; in particular, integer partitions. \\

\begin{figure}
	\includegraphics[scale=0.23]{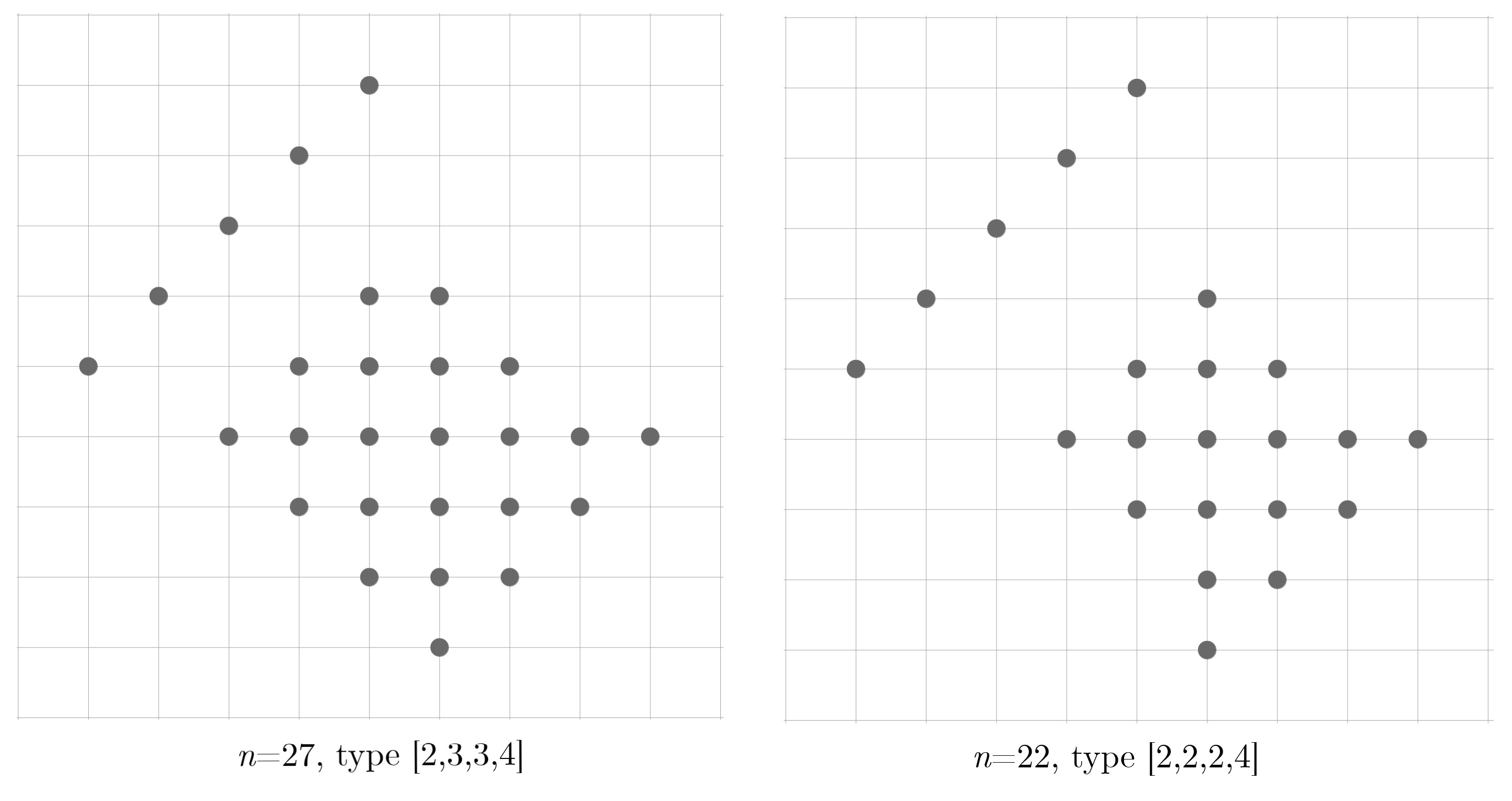} 
	\caption{The unique inseparable odd-sized set with $n\leq 41$ points (left)
		and the smallest even-sized set not of type $[N,M,N,M]$ (right; see Proposition \ref{prop-struc}).}
	\label{ger-p}
\end{figure}

Recall that the  symmetry group of a plane square lattice is $\mathbb{G}=\mathbb{D}_{4}\ltimes\mathbb{Z}^{2}$, denoted in crystallography as $p4m$. Since only monotone paths will be considered, the adjective ``monotone" will be skipped altogether. Paths will be denoted by curly letters $\mathscr{A},\mathscr{B}$, and so on. A finite lattice subset of size $n$ without a friendly path will be called an \emph{inseparable} $n$-set.
\begin{thmm}The following statements hold.
	\begin{itemize}
		\item[i)] Suppose $n\in\mathbb{N}$ is odd, $n\leq 41$, and a set of $n$ lattice points is inseparable. Then $n=27$. Up to the $\mathbb{G}$-symmetry, the unique configuration is presented in Figure \ref{ger-p} (left).
		\item[ii)] An inseparable $n$-set exists for every even number $n\geq 12$.
		\item[iii)]Odd numbers $n$ for which an inseparable $n$-set exists have a natural density $1$ among all odd numbers.
	\end{itemize}
\end{thmm}
\begin{defin}Let $c(n)$ and $\hat{c}(n)$ be, up to a $\mathbb{G}$-symmetry, and, respectively, $\mathbb{Z}^{2}$-symmetry, a number of different inseparable $n$-sets.  
\end{defin}
We will soon ascertain that this number is indeed finite.  These two sequences are related via an identity
\begin{eqnarray*}
	\hat{c}(n)=\sum\limits_{{\text{Over all  }\mathcal{P},\,|\mathcal{P}|=n,}\atop\text{Up to the $\mathbb{G}$-symmetry}}\frac{8}{g(\mathcal{P})},
\end{eqnarray*}
where $g(\mathcal{P})$ is the number of symmetries of a set $\mathcal{P}$. Obviously, $g(\mathcal{P})=1,2,4$, or $8$. \emph{A posteriori} it appears that even-sized and odd-sized cases have no conceptual difference, save when dealing with symmetry. The construction carried out in the next subsection (namely, a quartering) will immediately show that
\begin{eqnarray*}  
	\hat{c}(n)=8c(n)\text{ for odd }n.
\end{eqnarray*}	
\indent For small $n$, symmetric inseparable sets dominate and thus only even-sized sets exist. As $n$ grows larger, the portion of symmetric sets diminish and both cases tend to exhibit more and more complex configurational variety. In particular, $n=44$ is the first even number possessing an $n$-set with $g(\mathcal{P})=1$ (see Figure \ref{var4}, right).\\

\begin{table}
	\begin{tabular}{|r|r|r || r | r |  r || r| r| r|| r| r| r|}
		\hline
		\multicolumn{12}{|c|}{$\displaystyle{n\,|\,c(n)\,|\,\hat{c}(n)\,}$}\\		
		\hline
		
		\hline
			\quad	$1$ & $0$  & $\mathbf{0}$ &\quad\quad  $16$ &  $1$ & $\mathbf{2}$ &\quad\quad  $31$ & $0$ & $\mathbf{0}$ &\quad\quad    $46$ & $9$ & $\mathbf{30}$\\
		\hline

		\quad	 $2$&  $0$  & $ \mathbf{0}$ &\quad\quad   $17$&  $0$ & $\mathbf{0}$ &\quad\quad  $32$ & $3$ & $\mathbf{8}$ &\quad\quad    $47$ & $0$ & $\mathbf{0}$\\
		\hline
		
		\quad	 $3$ &  $0$ & $\mathbf{0}$ &\quad\quad  $18$ &  $2$ & $\mathbf{4}$ &\quad\quad  $33$ & $0$ & $\mathbf{0}$ &\quad\quad    $48$ & $5$ & $\mathbf{13}$\\
			\hline
		
		\quad	 $4$ &  $1$ & $\mathbf{0}$ &\quad\quad   $19$ &  $0$ & $\mathbf{0}$ &\quad\quad  $34$  & $6$ & $\mathbf{18}$ &\quad\quad   $49$  & $0$ & $\mathbf{0}$\\
			\hline
		
		\quad	 $5$ &  $0$ & $\mathbf{0}$ &\quad\quad   $20$ &  $2$ & $\mathbf{6}$ &\quad\quad  $35$ & $0$ & $\mathbf{0}$ &\quad\quad    $50$ & $7$ & $\mathbf{22}$\\
			\hline
		
		\quad	 $6$ &  $0$ & $\mathbf{0}$ &\quad\quad  $21$ &  $0$ & $\mathbf{0}$ &\quad\quad  $36$ & $10$ & $\mathbf{32}$ &\quad\quad    $51$ & $0$ & $\mathbf{0}$\\
			\hline
		
		\quad	 $7$ &  $0$ & $\mathbf{0}$ &\quad\quad  $22$ &  $2$ & $\mathbf{6}$ &\quad\quad  $37$ & $0$ & $\mathbf{0}$ &	\quad\quad   $52$ & $8$ & $\mathbf{32}$\\
			\hline
		
		\quad	 $8$ &  $1$ & $\mathbf{2}$ &\quad\quad  $23$ &  $0$ & $\mathbf{0}$ &\quad\quad  $38$ & $6$ & $\mathbf{20}$ &\quad\quad   $53$  & $0$ & $\mathbf{0}$\\
			\hline
		
		\quad	 $9$ &  $0$ & $\mathbf{0}$ &\quad\quad  $24$ &  $4$ & $\mathbf{11}$ &\quad\quad  $39$ & $0$ & $\mathbf{0}$ &\quad\quad   $54$  & $12$ & $\mathbf{42}$\\
			\hline
		
		\quad	 $10$ &  $0$ & $\mathbf{0}$ &\quad\quad  $25$ &  $0$ & $\mathbf{0}$ &\quad\quad  $40$ & $9$ & $\mathbf{29}$ &\quad\quad   $55$  & $0$ & $\mathbf{0}$\\
			\hline
		
		\quad	 $11$ &  $0$ & $\mathbf{0}$ &\quad\quad  $26$ &  $4$ & $\mathbf{12}$ &\quad\quad  $41$ & $0$ & $\mathbf{0}$ &\quad\quad   $56$  & $18$ & $\mathbf{64}$\\
			\hline
		
		\quad	 $12$ &  $1$ & $\mathbf{1}$ &\quad\quad  $27$ &  $1$ & $\mathbf{8}$ &\quad\quad  $42$ & $12$ & $\mathbf{42}$ &\quad\quad   $57$  & $0$ & $\mathbf{0}$\\
			\hline
		
		\quad	 $13$ &  $0$ & $\mathbf{0}$ &\quad\quad  $28$ &  $3$ & $\mathbf{7}$ &\quad\quad  $43$ & $1$ & $\mathbf{8}$ &\quad\quad   $58$  & $14$ & $\mathbf{50}$\\
			\hline
		
		\quad	 $14$ &  $1$ & $\mathbf{2}$ &\quad\quad  $29$ &  $0$ & $\mathbf{0}$ &\quad\quad  $44$ & $18$ & $\mathbf{67}$ &\quad\quad   $59$  & $0$ & $\mathbf{0}$\\
			\hline
		
		\quad	 $15$ &  $0$ & $\mathbf{0}$ &\quad\quad  $30$ &  $2$ & $\mathbf{6}$ &\quad\quad  $45$ & $2$ & $\mathbf{16}$ &\quad\quad   $60$  & $17$ & $\mathbf{64}$\\ 
		
		\hline \hline
	\end{tabular}
	\caption{Sequences $c(n)$ and $\hat{c}(n)$ (in bold), $1\leq n\leq 60$.}
	\label{tab}
\end{table}

\begin{figure}
	\includegraphics[scale=0.17]{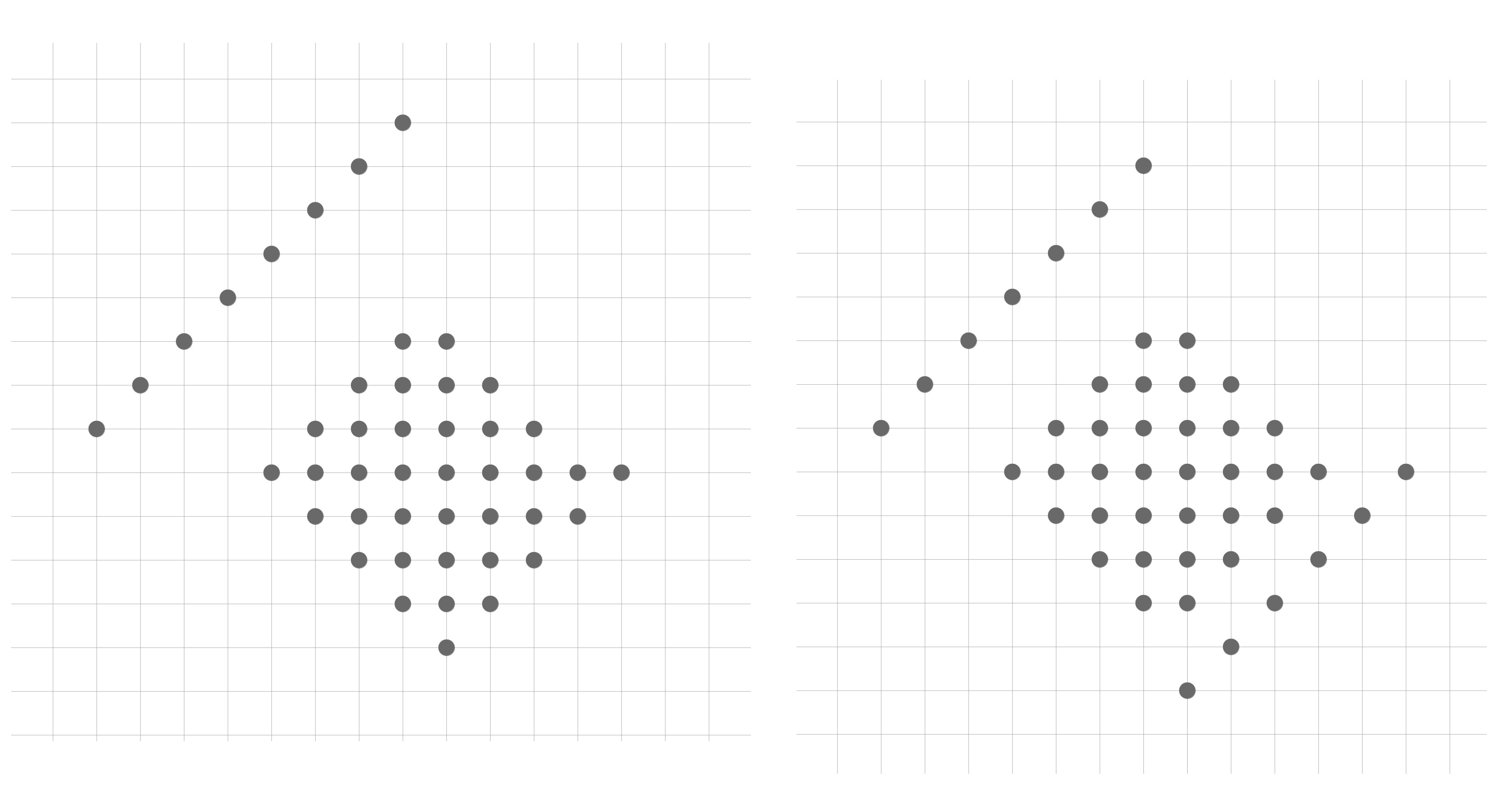} 
	\caption{Two inseparable examples for $n=45$}
	\label{ex45}
\end{figure}
   
Prior to passing to proofs, we note that the given setup has a natural $N$-dimensional analogue for all positive integers $N\geq 2$. Also, it has a $2$-dimensional analogue for a hexagonal lattice. Though a topic of lattice paths is a basic one in combinatorics and is essentially covered in \cite{stanley}, the interaction between finite subsets of $\mathbb{Z}^{2}$ and lattice paths, provided by the notion of a ``friendly path", is a new and intriguing topic with several ramifications (possibly, even leading to modular forms; see Subsection \ref{expect}).  
   
\subsection{Preliminary quartering}
For a path $\mathscr{A}$, define \emph{the balance} of $\mathscr{A}$ by $d(\mathscr{A})=\ell(\mathscr{A})-r(\mathscr{A})$, where $\ell$ and $r$ are number of points from $\mathcal{P}$ on the left and, respectively, on the right shores of $\mathscr{A}$. Let $T(\mathscr{A})$ be the number of points on the path itself (this will be needed in Section \ref{sec-fin}). A \emph{turn} on the path can be either left or right. These will be called \emph{stop-points}. If $\mathscr{C}$ and $\mathscr{D}$ are two paths such that all stop-points of 
$\mathscr{C}$ lie on the left shore of $\mathscr{D}$ or $\mathscr{D}$ itself, we write $\mathscr{C}\precsim\mathscr{D}$. This implies $d(\mathscr{C})\leq d(\mathscr{D})$. The symbol $``\precsim"$ thus gives a partial ordering of paths.  Let $\mathbf{S}$ be a closed rectangle in the plane bounded by lines $X=x_{0}$, $X=x_{1}$, $Y=y_{0}$, $Y=y_{1}$, $x_{0}<x_{1}$, $y_{0}<y_{1}$, such that $\mathcal{P}\subset\mathbf{S}$. We write $\mathscr{A}\doteq\mathscr{B}$, if $\mathscr{A}\cap\mathbf{S}=\mathscr{B}\cap\mathbf{S}$. Naturally, the structure of the path outside $\mathbf{S}$ does not affect its relation towards $\mathcal{P}$.\\

Let $n\in\mathbb{N}$. Suppose that a set $\mathcal{P}$, $|\mathcal{P}|=n$, is inseparable. First, we will show that its structure can be fully described with a help of an identity $n=\sum_{j=1}^{4}Y_{j}$, $Y_{j}\in\mathbb{N}$, and four sets $L_{j}$, $1\leq j\leq 4$. Each $L_{j}$ is a certain partition of $Y_{j}$ into distinct positive integral parts (see A000009 in \cite{oeis}). More details will be given by Proposition \ref{prop-struc}.\\

Consider a horizontal path $y=K$, $K\in\mathbb{Z}$. Horizontal paths (going left to right) are both uphill and downhill. If $K$ is large enough, a path has a balance $-n$. Similarly, if $L$ is large enough, a horizontal line $y=-L$,  $L\in\mathbb{Z}$, has a balance $n$. Since the sequence $d\big{(}[y=Y]\big{)}$ for $Y$ ranging from $K$ to $-L$ is non-decreasing, there exists a place where this sequence changes sign. By our assumption no path has a balance $0$. This implies the existence of $M\in\mathbb{Z}$ such that
\begin{eqnarray}
d\Big{(}[y=M+1]\Big{)}<0,\quad d\Big{(}[y=M]\Big{)}>0.\label{red-hor}
\end{eqnarray}
Let us draw a red line $y=M+\frac{1}{2}$. In the same manner, consider vertical lines going upwards. They are uphill. We thus find $T\in\mathbb{Z}$ such that
\begin{eqnarray}
d\Big{(}[x=T]\Big{)}<0,\quad d\Big{(}[x=T+1]\Big{)}>0.\label{red-vert}
\end{eqnarray}

\begin{figure}
\includegraphics[scale=.47]{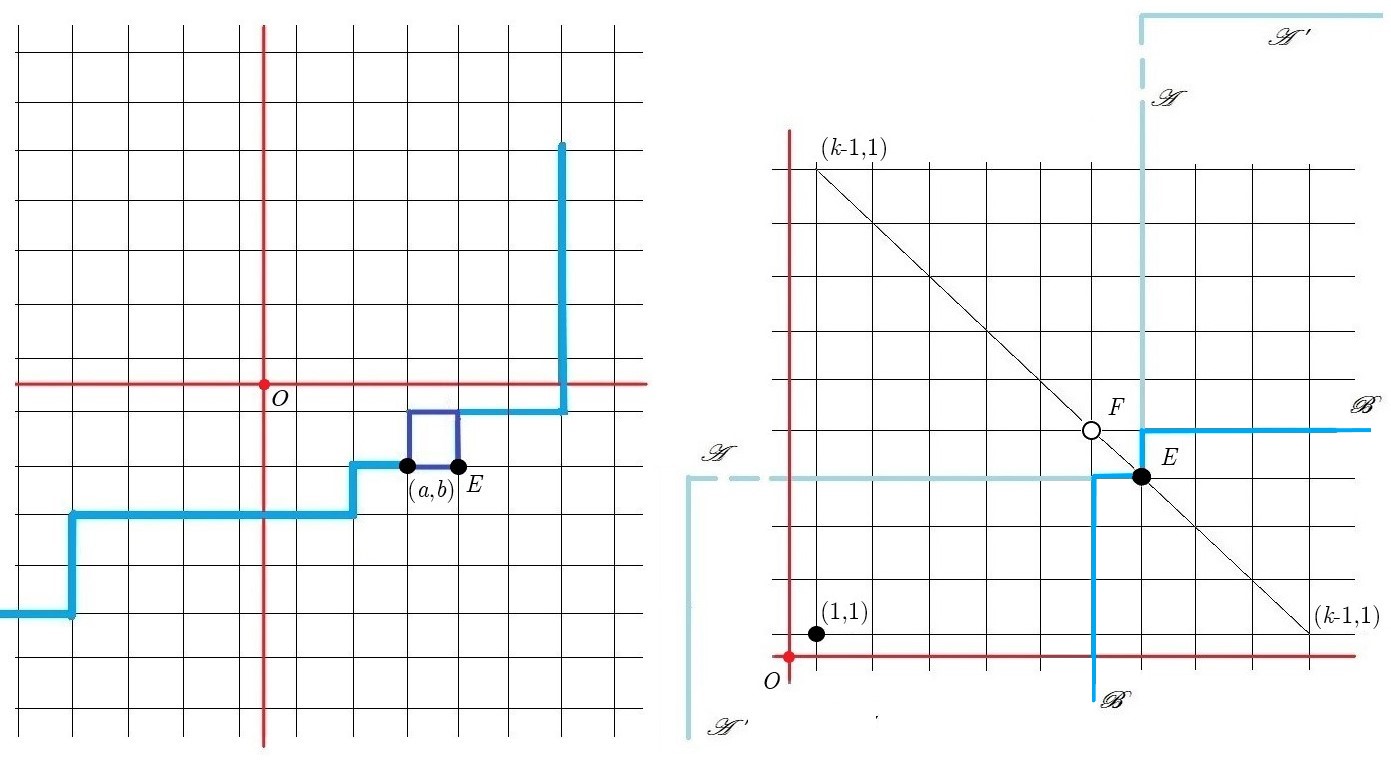} 
\caption{Quartering of the set $\mathcal{P}$}
\label{quart}
\end{figure}

\begin{figure}
\includegraphics[scale=0.28]{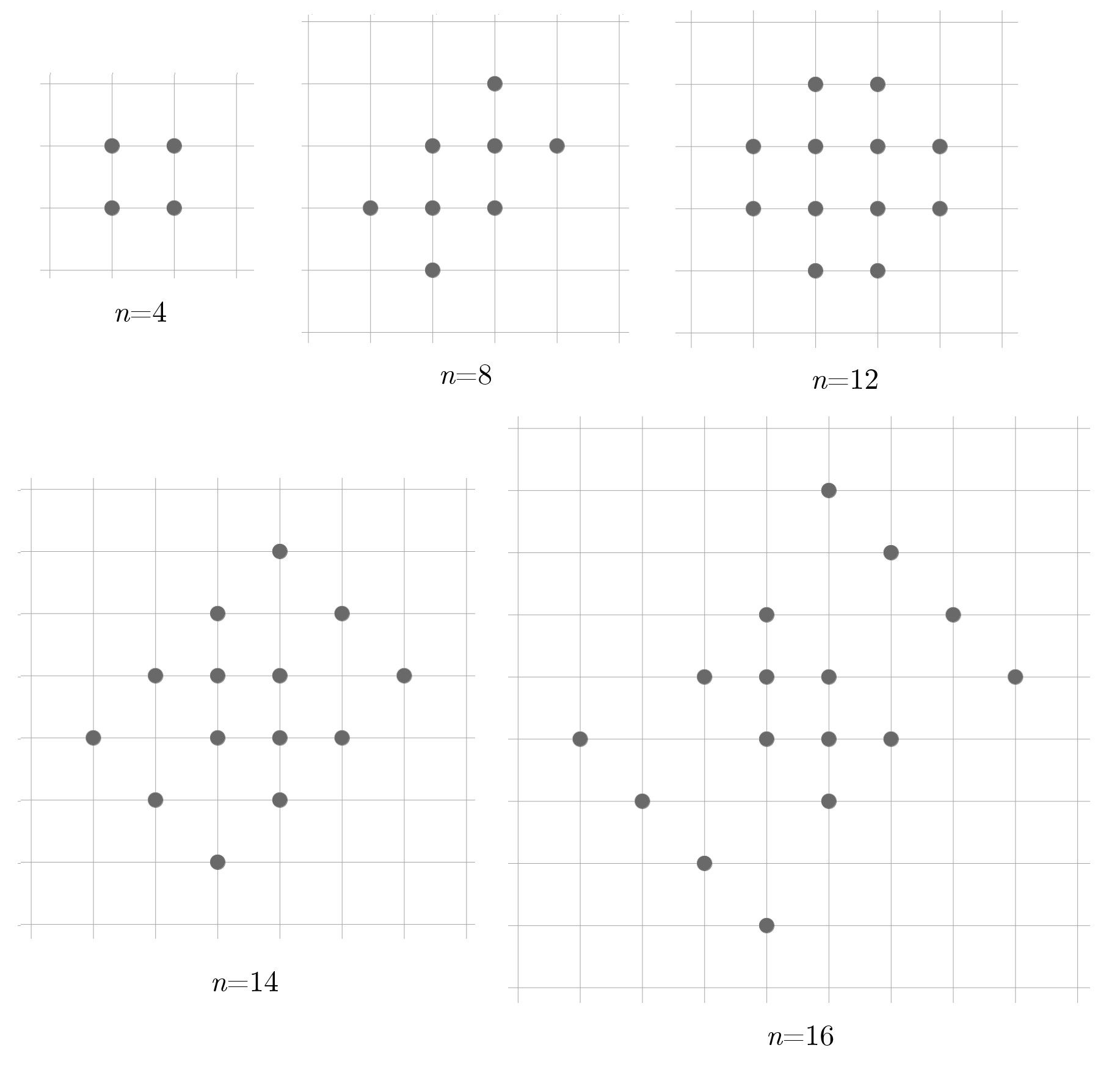} 
\caption{Examples of all even $n$ with $c(n)=1$}
\label{ex1}
\end{figure}

\begin{figure}
\includegraphics[scale=0.20]{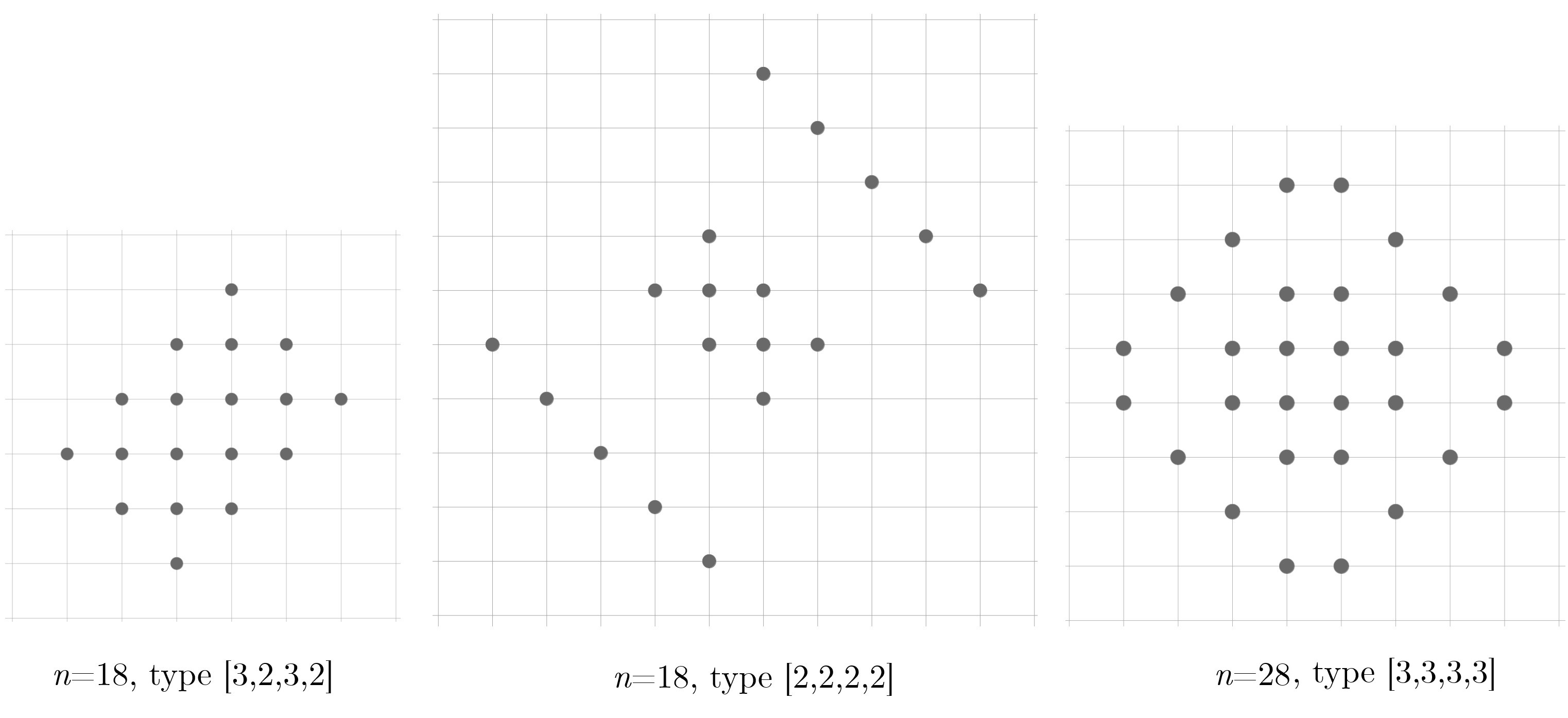} 
\caption{Two inseparable examples for $n=18$ and the the smallest lacunary inseparable configuration with $8$-fold symmetry ($n=28$)}
\label{ex2}
\end{figure}

\begin{figure}
\includegraphics[scale=0.18]{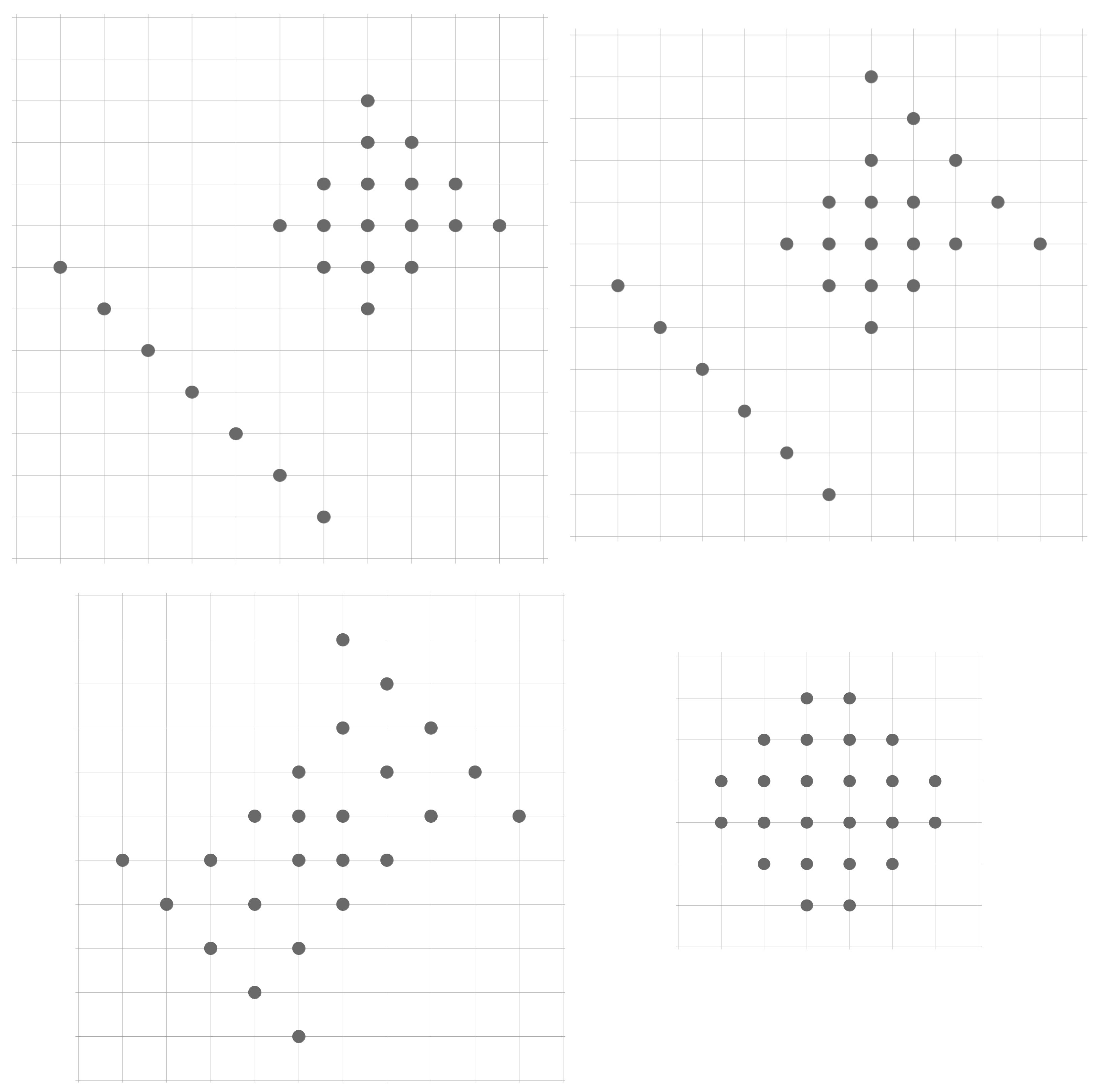} 
	\caption{Four inseparable examples for $n=24$}
\label{ex3}
\end{figure}

Let us draw a vertical red line $x=T+\frac{1}{2}$. Let $O$ be the intersection of these two. This unambiguous construction quarters points in $\mathcal{P}$.\footnote{This method also provides an algorithm for finding a friendly path. If vertical or horizontal line can be such, we are done. In not, there exists the unique quartering of points in $\mathcal{P}$. The second step of the algorithm becomes clear minding the analysis in Section \ref{analit}.} The group of isometries of the plane which leaves the union of two red lines intact is $\mathbb{D}_{4}$. If $\sigma\in\mathbb{D}_{4}$, then trivially $\sigma(\mathcal{P})$ is inseparable, too.

\section{Representations of inseparable sets}
\label{analit}
 \subsection{The structure}We now define a procedure on the path called \emph{a shift}, which can be \emph{negative} or \emph{positive}. Take any left-turn of the uphill path. Two unit segments at this stop-point are $(a,b)\rightarrow(a+1,b)$ and $(a+1,b)\rightarrow(a+1,b+1)$. Let us replace them by $(a,b)\rightarrow(a,b+1)$ and $(a,b+1)\rightarrow(a+1,b+1)$ (see Figure \ref{quart}, left, left-turn at $E$). After this replacement the balance of the path changes by $0$, $-1$, or $-2$. This  is a \emph{negative shift}. Let us also define a negative shift for a downhill path. In this case, a left turn need also be taken, and two segments $(a,b)\rightarrow(a,b-1)$ and $(a,b-1)\rightarrow(a+1,b-1)$ are replaced by $(a,b)\rightarrow(a+1,b)$ and $(a+1,b)\rightarrow(a+1,b-1)$. This changes the balance also by $0$, $-1$, or $-2$. On the other hand, if we consider right-turns, we can analogously introduce the notion of a \emph{positive shift}. The latter changes the balance by $0$, $1$, or $2$.\\

\indent Let us return to the quartering of the set $\mathcal{P}$. Consider the first (top-right) quadrant $\mathcal{Q}_{1}$. Let the lattice points there have coordinates $(x,y)$, $x,y\in\mathbb{N}$ (see Figure \ref{quart}, right). For fixed $k\in\mathbb{N}$, $k\geq 1$, take an integer diagonal $\{(x,k+1-x): 1\leq x\leq k\}$. Call this set $\Upsilon_{1,k}$. We will now prove the following statement.
\begin{lemm}
All or none points of this diagonal belong to $\mathcal{P}$.
\end{lemm}
For $k=1$ this is a tautology. Suppose there exists $k\geq 2$ and two neighbouring points $E,F$ on this diagonal satisfying $E\in\mathcal{P}$, $F\notin\mathcal{P}$. Without loss of generality, assume $E$ is to the right of $F$. Otherwise we just use the reflection with respect to the line $x=y$. So, $E=(x,y)$, $F=(x-1,y+1)$, $x\geq 2$.\\

\indent Consider now two uphill paths $\mathscr{A}$ and $\mathscr{B}$, as marked in Figure \ref{quart}, right. Both dotted half-lines extend to infinity. Outside $\mathbf{S}$, the path $\mathscr{A}$ can be altered to attain a shape of a path $\mathscr{A}'$. This does not change the balance. Now, inside $\mathbf{S}$ the path $\mathscr{A}'$ can be gradually transformed into the the path $\mathscr{B}$ by a series of positive shifts, leaving the left-turn at $E$ intact:
\begin{eqnarray*}
\mathscr{A}'=\mathscr{E}_{0}\precsim\mathscr{E}_{1}\precsim\cdots\mathscr{E}_{r}\doteq\mathscr{B}.
\end{eqnarray*}
Moreover, since $\mathscr{A}\precsim[y=M+1]$, we have $d(\mathscr{A}')=d(\mathscr{A})<0$. In a similar vein, since $[x=T+1]\precsim\mathscr{B}$, $d(\mathscr{B})>0$ (here we implicitly use a condition $x\geq 2$). By our assumption, a sequence  $\{d(\mathscr{E}_{i}): 0\leq i\leq r\}$ is monotone, does not contain $0$ and makes integer jumps by at most $2$. Since it goes from negative to a positive value, there exists $s$, $1\leq s\leq r$, such that $d(\mathscr{E}_{s-1})=-1$ and $d(\mathscr{E}_{s})=1$. Consider the path $\mathscr{E}_{s}$. Making a negative shift for $\mathscr{E}_{s}$  at $E$ produces a path with a balance $0$. A contradiction. $\blacksquare$\\

The same conclusion must hold for every integer diagonal of any other quadrant. Thus, we have proved the main structural result. 
\begin{prop}
Every inseparable set $\mathcal{P}$ can be represented in the form 
\begin{eqnarray}
\theta=[L_{1}L_{2}L_{3}L_{4}],&& L_{r}=(v_{r,1},v_{r,2},\ldots,v_{r,N_{r}}),\label{repr}\\
1\leq v_{r,1}<v_{r,2}<\cdots <v_{r,N_{r}},&& r\in\{1,2,3,4\},\,\, N_{r}\in\mathbb{N}.\nonumber
\end{eqnarray} 
Here $1,2,3$ and $4$ correspond, respectively, to the namesake quadrant. The notation implies that integer diagonals with numbers $v_{r,i}$ fully belong to $\mathcal{P}$, and these exhaust all of its points. We have
\begin{eqnarray*}
n=n(\theta)=\sum\limits_{r=1}^{4}Y_{r},\text{ where } Y_{r}=\sum\limits_{t=1}^{N_{r}}v_{r,t}.
\end{eqnarray*}
\label{prop-struc}
\end{prop}
Let us call the quadruple $[N_{1},N_{2},N_{3},N_{4}]$, the \emph{type} of the representation $\theta$, and $\max N_{i}-\min N_{i}$ its \emph{variation}.\\

Calculations will confirm that in fact $v_{r,1}=1$ for all inseparable sets. Surely, only specific choices of such representation of an integer $n$ produces an inseparable configuration. Computational part of this proposition will be treated in Section \ref{sec-fin}, while a theoretical part will be dealt with in \cite{antra-dalis}.  We however remark that all results of this paper can be checked by hand and do not depend on computer codes. 
\subsection{Examples}
\begin{figure}
	\includegraphics[scale=0.24]{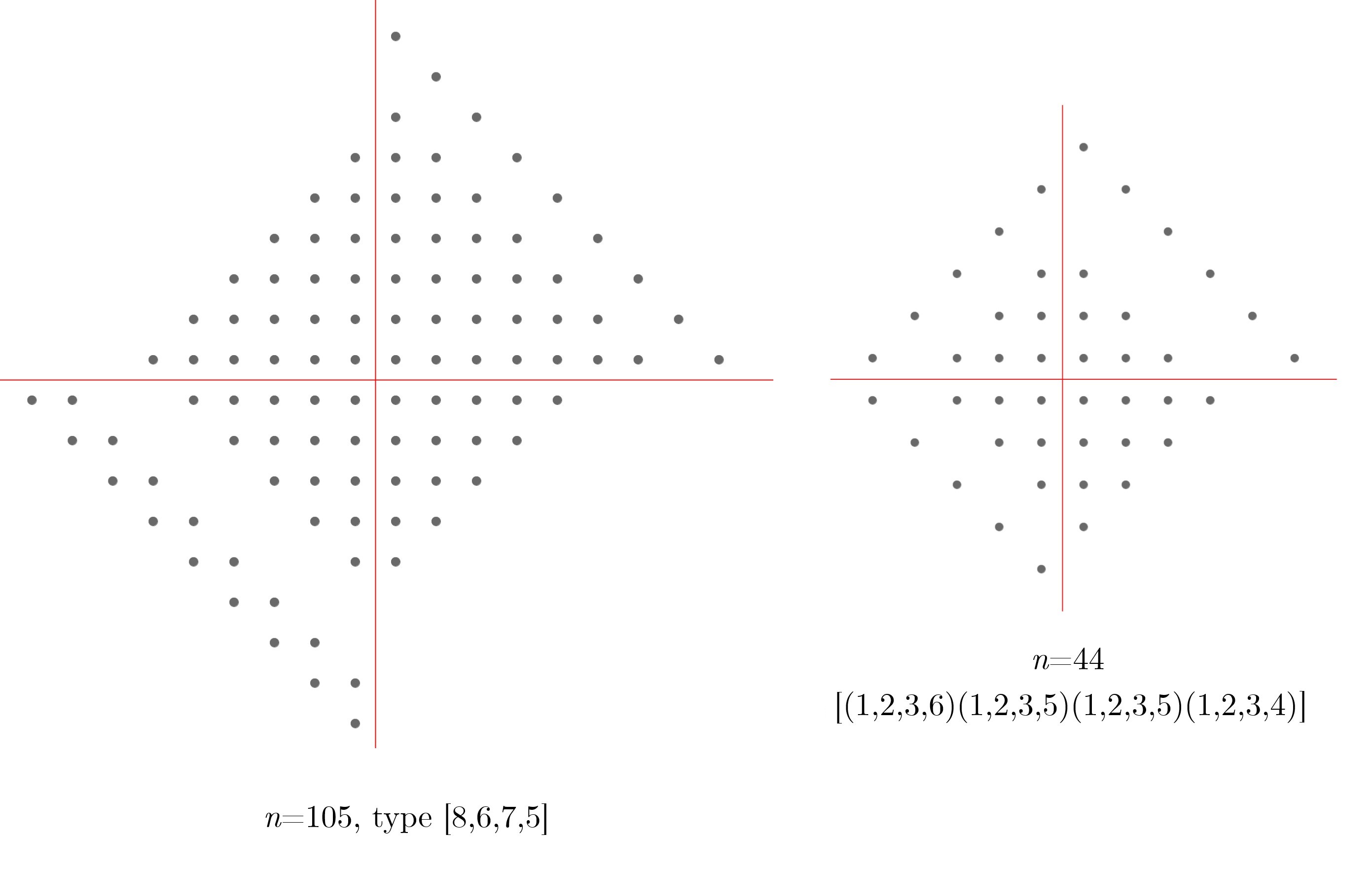}
	\caption{The smallest odd-sized inseparable set with variation $3$ (left),\\
		and the smallest even-sized inseparable set without symmetries (right).}
	\label{var4}
\end{figure}

The smallest integer with $c(n)\geq 2$ is $n=18$ (see Figure \ref{ex2}), where these collections are, respectively,
\begin{eqnarray*}
\theta_{1}=[(1,2,3)(1,2)(1,2,3)(1,2)]\text{ and }\theta_{2}=[(1,5)(1,2)(1,5)(1,2)].
\end{eqnarray*} 
The smallest integer with $c(n)\geq 3$ is $n=24$ (Figure \ref{ex3}), where in fact $c(24)=4$:
\begin{center}
	\begin{tabular}{ l l }
$\theta_{1}=[(1,3,5)(1,2)(1,3,5)(1,2)]$,& $g(\theta_{1})=4$;\\
$\theta_{2}=[(1,2,3,5)(1,2)(1,6)(1,2)]$,& $g(\theta_{2})=2$;\\
$\theta_{3}=[(1,2,3,4)(1,2)(1,7)(1,2)]$,& $g(\theta_{3})=2$;\\
$\theta_{4}=[(1,2,3)(1,2,3)(1,2,3)(1,2,3)]$,& $g(\theta_{4})=8$.
	\end{tabular}
\end{center}
Consequently, 
\begin{eqnarray*}
	\hat{c}(24)=\frac{8}{4}+\frac{8}{2}+\frac{8}{2}+\frac{8}{8}=11.
	\end{eqnarray*}
The smallest odd-sized set with a variation $3$ is presented in Figure \ref{var4} (left). The initial solution in {\sc Monthly} shows that variation for even-sized inseparable sets can be arbitrarily large. With all tools at hand, we easily see that the same holds for odd-sized sets.

\section{Full conditions for inseparability}
\label{sec-fin}
Assume, we run a computer program which enumerates all possible quadruples $[L_{1}L_{2}L_{3}L_{4}]$ of finite increasing integer sequences. Which of them truly represent an inseparable set? Three tests and a preliminary requirement should be passed. This is the complete list of conditions.
\subsection{Evenness and proper quartering} 
\begin{center}
	\setlength{\shadowsize}{2pt}\shadowbox{\texttt{Test 0}}
\end{center}
We start from type $[N_{1},N_{2},N_{3},N_{4}]$ and $n(\theta)$ which (in principal) can lead to an inseparable set. By this it is meant that the following condition holds: 
\begin{eqnarray*}
	N_{1}+N_{3}\equiv N_{2}+N_{4}\equiv n\text{ (mod }2).
\end{eqnarray*}
If satisfied, this implies that the generic path starting in quarter $\mathcal{Q}_{3}$ and ending in quarter $\mathcal{Q}_{1}$ (the same applies to the pair $\mathcal{Q}_{2}$ and $\mathcal{Q}_{4}$) pass through $N_{1}+N_{3}+1$ points. Therefore, the remaining $n-N_{1}-N_{3}-1$ points (an odd number) cannot be divided exactly in half. A friendly path, if it exists, must not be a generic one (see \textbf{Case 3} below for a precise meaning of this).\\
\begin{center}
	\setlength{\shadowsize}{2pt}\shadowbox{\texttt{Test A}}
	\end{center}
Next, we must check that the quartering of the set is proper. This means
\begin{itemize}
\item[\textbf{1) }]	$-N_{3}-N_{4}< Y_{1}+Y_{2}-Y_{3}-Y_{4}<N_{1}+N_{2}$;\\\
\item[\textbf{2) }]$	-N_{2}-N_{3}< Y_{1}+Y_{4}-Y_{2}-Y_{3}<N_{1}+N_{4}$.\\
\end{itemize} 
Using reflections with respect to both red lines in case of necessity, one can always assume that $Y_{1}+Y_{2}\geq Y_{3}+Y_{4}$ and $Y_{1}+Y_{4}\geq Y_{2}+Y_{3}$ holds. However, for computational purposes it is more convenient to work with conditions \textbf{1)} and \textbf{2)} directly. Depending on symmetries of a representation type $[N_{1},N_{2},N_{3},N_{4}]$, additional conditions are imposed by a computer program which checks all configurations of a  particular type for inseparability.\\

To better understand where the next two tests come from, consider an example. 
\subsection{Example with $n=108$}
\label{simtast}
\begin{figure}
	\includegraphics[scale=0.16]{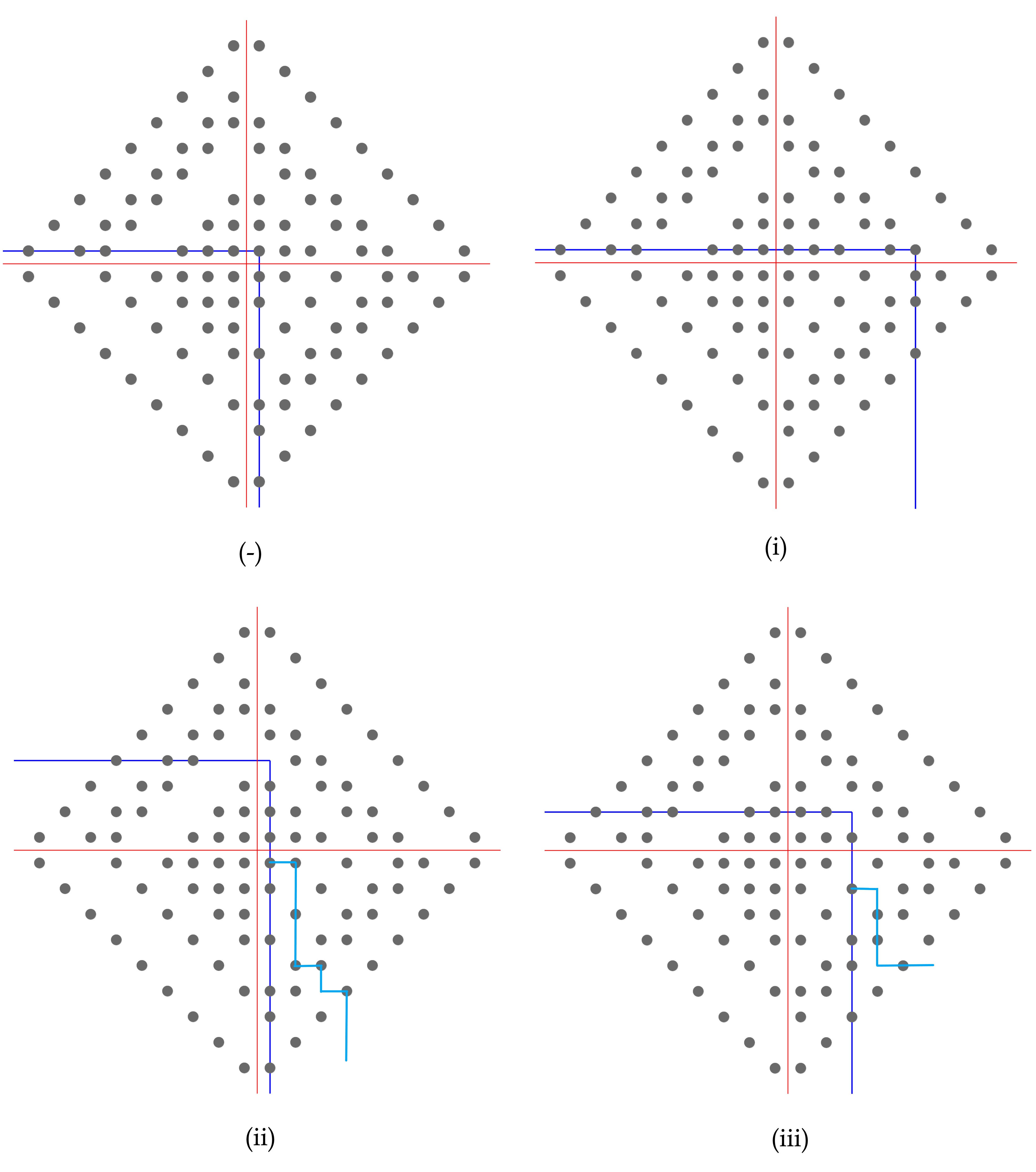} 
	\caption{An example with $n=108$}
	\label{ex108}
\end{figure}

Let $\theta=[L_{1}L_{2}L_{3}L_{4}]$, where
\begin{eqnarray*}
L_{1}=(1,2,3,5,6,9),\quad L_{2}=(1,2,3,6,7,9),\quad L_{3}:=(1,2,3,4,6,9),\quad 
L_{4}:=(1,2,4,6,7,9).
\end{eqnarray*}
Here  $Y=(26,28,25,29)$, $(N_{1},N_{2},N_{3},N_{4})=(6,6,6,6)$; see Figure \ref{ex108}. \\

`\textbf{Case 3}'. For this particular set, assume there exists a friendly path which enters the configuration at $\mathcal{Q}_{2}$, passes through $\mathcal{Q}_{1}$ and exists it at $\mathcal{Q}_{4}$. So, this friendly path have no common points with $\mathcal{Q}_{3}$. This is exactly what the name `\textbf{Case 3}' is supposed to mean (`Case 1', `Case 2', `Case 4' are defined analogously). Now come several subcases.\\

\textbf{`Subcase 3-'}. Consider a downhill path which takes a single right-turn at the point $(1;1,1)$ (the first coordinate of this notation indicates the index of the quarter; see Figure \ref{ex108}, top-left). There are $N_{2}+N_{4}+1=13$ points on this path. Consequently, it cannot be friendly.\\

\textbf{`Subcase 3i'}. We now ask the following: find the smallest $x_{0}\in\mathbb{N}$ such that a downhill path, which makes a single right-turn at $(1;x_{0},1)$, contains an even number of points (this is what ``non generic" means). Only such paths can potentially be friendly. The answer is $x_{0}=6$ (Figure \ref{ex108}, top-right). For this particular path, $T(\mathscr{A})=14$, $r(\mathscr{A})=47$, $\ell(\mathscr{A})=47$. By a lucky chance, it is friendly.\\

\textbf{`Subcase 3ii'}. In the same vein, if a downhill path makes a single right-turn at $(1;1,y_{0})$, it contains an even amount of points, and $y_{0}$ is the smallest possible, then $y_{0}=4$ (same picture, bottom-left).  For this path (shown in blue), $T(\mathscr{A})=12$, $r(\mathscr{A})=40$, $\ell(\mathscr{A})=56$.\\

\textbf{`Subcase 3iii'}. As a final sub-case, note that the first positive integer $r$ not among members of $L_{1}$ is $r=4$. Consequently, any down-hill path making a single right-turn at points $(1;3,2)$ or $(1;2,3)$ contains and even amount points ($12$ each). In both cases, $r(\mathscr{A})=44$, $\ell(\mathscr{A})=52$ (same picture, bottom-right). The second path is shown in blue.\\

Now, a negative shift made strictly inside $\mathcal{Q}_{2}$ or $\mathcal{Q}_{4}$ alters a balance by $0$ or $-2$ ($1$ is excluded due to construction). Therefore, in the last two subcases we can construct a friendly path starting with a blue path by a series of negative shifts. Results in both cases are shown in sky-blue. For small $n$ it might happen that, while performing a series of negative shifts, we are forced to leave $\mathcal{Q}_{2}$ or $\mathcal{Q}_{4}$. Yet a distinguished right turn chosen controls a positive shift by $+1$. Thus, a friendly path is always to be found.\\

This gives a clear picture how one checks for inseparability in a general case: in \textbf{Case 3} and all sub-cases (the path is downhill), inequalities $r(\mathscr{A})>\ell(\mathscr{A})$ must hold.\\

`\textbf{Subcase 3iv}'. A special attention needs to be given if $v_{r,1}>1$ for some $r$. In this case \texttt{test B} must be adjusted. Suppose, $v_{1,1}>1$ holds (Figure \ref{fails}, right). Since $10=r(\mathscr{A})=Y_{3}>Y_{1}+Y_{2}+Y_{4}-N_{2}-N_{4}=\ell(\mathscr{A})=8$, this does not lead to a friendly path. Assume now that instead we had an inequality $r(\mathscr{A})<\ell(\mathscr{A})$ . The tricky part is that, differently from `\textbf{Subcase 3ii}', we cannot guarantee that with a help of a series of negative shifts this path can be transformed into a friendly one. This negative outcome occurs precisely if the balance of a sky-blue path is also positive: $Y_{1}>Y_{2}+Y_{3}+Y_{4}-N_{2}-N_{4}$. Thus, in case $v_{1,1}>1$ a friendly path passing through $(1;1,1)$ does not exist if and only if
\begin{eqnarray}
	|Y_{3}-Y_{1}|>Y_{2}+Y_{4}-N_{2}-N_{4}.\label{iskr}
\end{eqnarray} 
In case of (Figure \ref{fails}, right), a friendly path does exist (the one which avoids $\mathcal{Q}_{2}$), but the condition $v_{1,1}>1$ cannot be ruled out for inseparable sets with small $n$. Nevertheless, computations confirm that indeed no such set exists. 
\subsection{General check} Consider general $n$ and a subcase \textbf{3ii}. How we find such $x_{0}$? Checking many particular examples, we arrive at the following conclusion. For each ordered pair $(\mathcal{Q}_{i},\mathcal{Q}_{j})$ of neighbouring quadrants (that is, $i-j\equiv \pm 1\text{ (mod }4)$), let $w_{i,j}$ be the smallest $p\in\mathbb{N}\setminus\{1\}$ such that 
 \begin{eqnarray*}
\text{either }p\in L_{i},\, p-1\notin L_{j},\text{ or }p\notin L_{i},\, p-1\in L_{j}. 
 \end{eqnarray*}
Note that such $p$ might not exist. In this case for all paths $\ell$ not passing through a quadrant $j+1$, $T(\ell)$ is odd. Otherwise, if such $p$ exists, this is exactly the quantity we are looking for (Looking back at the subcase \textbf{3ii} in the previous subsection, $w_{1,4}=6$ corresponds exactly to the point $x_{0}=6$). As we have just witnessed, if the set is inseparable, the inequality $\ell(\mathscr{A})<r(\mathscr{A})$ must hold for all pairs $(i,j)=(1,4), (4,3), (3,2), (2,1)$, and $r(\mathscr{A})<\ell(\mathscr{A})$ must hold for all pairs $(i,j)=(4,1), (1,2), (2,3), (3,4)$. MAPLE program which checks these criteria is given in the appendix. This is \texttt{Test B}. Its input is representation $\theta$ and a pair of integers $(j,i)$. The outcome is $\mathtt{[Ls,Rs,yra,Tt]}$. Here $\mathtt{Ls}=\ell(\mathscr{A})$, $\mathtt{Rs}=r(\mathscr{A})$, $\texttt{Tt}=T(\mathscr{A})$. Boolean variable $\texttt{yra}$ gives value $\texttt{false}$ if no such $p$ exists. In case $v_{r,1}>1$ this tests only checks the condition (\ref{iskr}). Naturally, the boolean output of this test is sufficient in all cases. We track values $\mathtt{Ls}$, $\mathtt{Rs}$ and $\mathtt{Tt}$ only for debugging purpose.    \\

Equally, analogue of $r$ in the subcase \textbf{3iii} is the following number. For quadrant $\mathcal{Q}_{i}$, let $r_{i}$ be the smallest $p\in\mathbb{N}$ not contained in $L_{i}$. MAPLE codes which checks that no such path can be deformed to become friendly is given in the Appendix as \texttt{Test C}. Figure \ref{fails} shows a configuration which passes \texttt{Test 0}, \texttt{Test A} and \texttt{Test B}, but fails \texttt{Test C}.\\
\begin{figure}
	\includegraphics[scale=0.24]{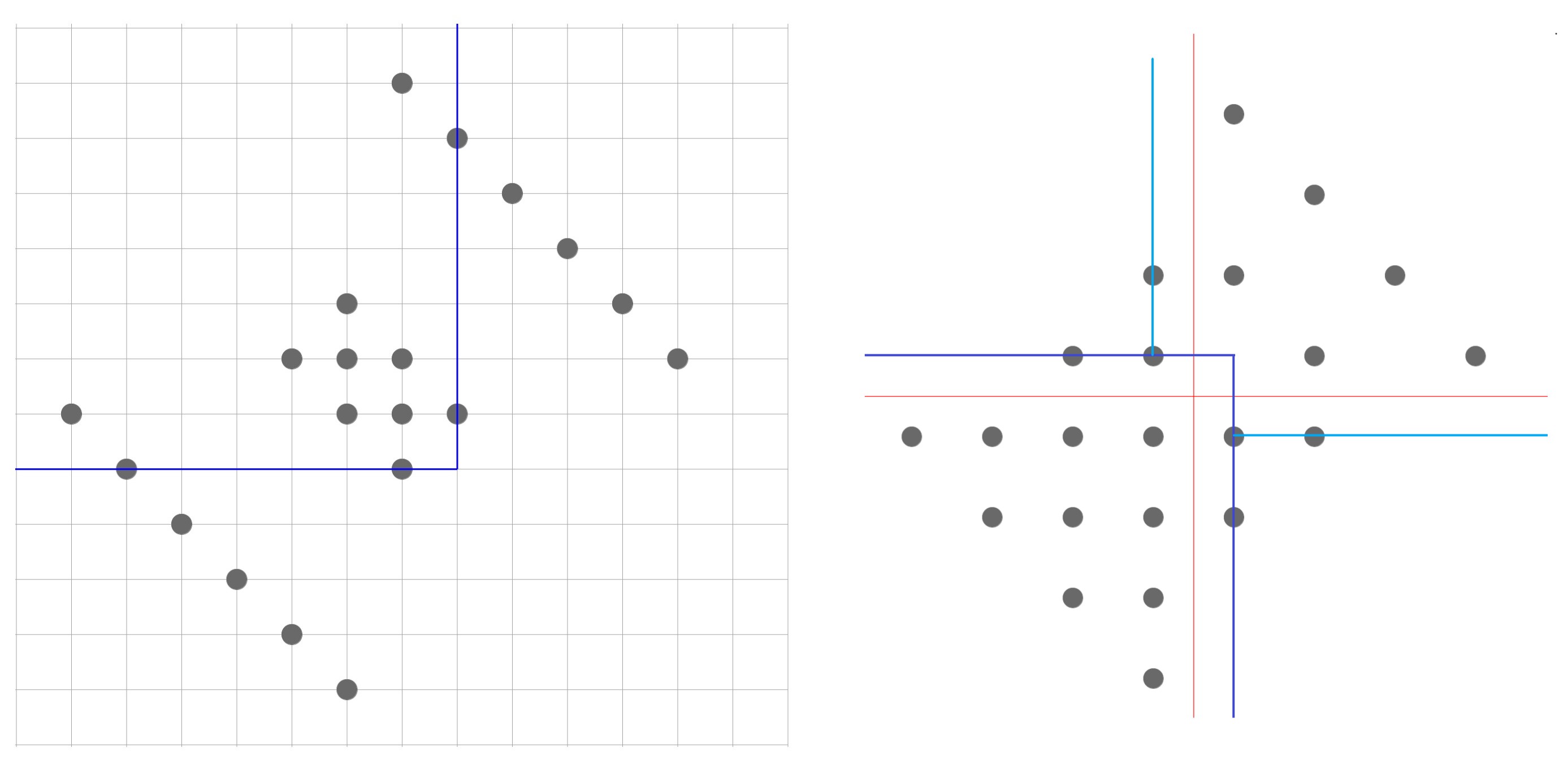} 
	\caption{$[(1,6)(1,2)(1,6)(1,2)]$ fails \texttt{Test C} (left) and supplementary illustration for \texttt{Test B} in case $v_{1,1}>1$ (right)}
	\label{fails}
\end{figure}

\section{Basic properties of both sequences}
\subsection{Even numbers with $c(n)=0$}
\begin{figure}
	\includegraphics[scale=0.16]{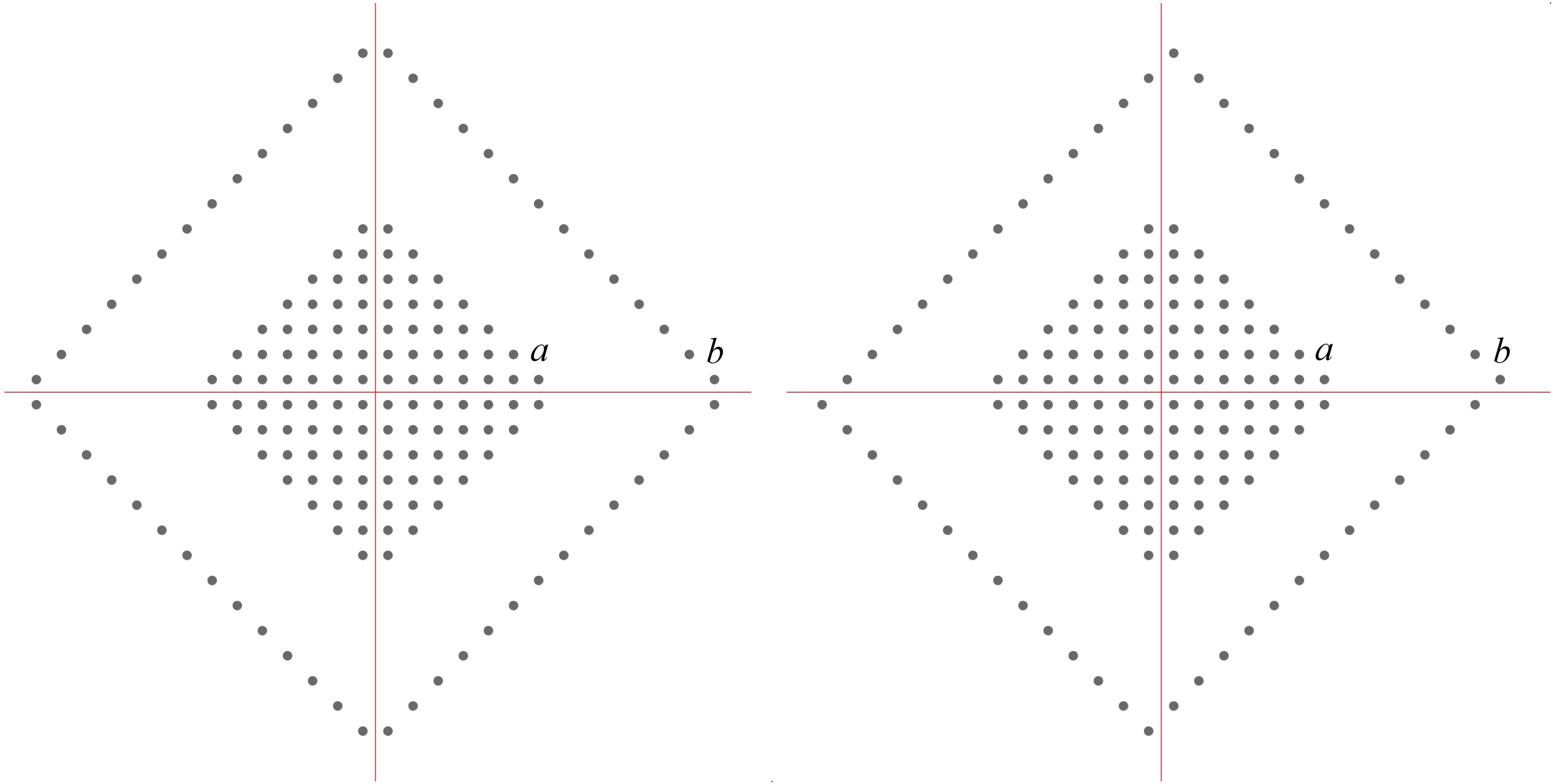} 
	\caption{$a+1\leq b\leq 2a$}
	\label{induct-1}
\end{figure}

\begin{figure}
	\includegraphics[scale=0.16]{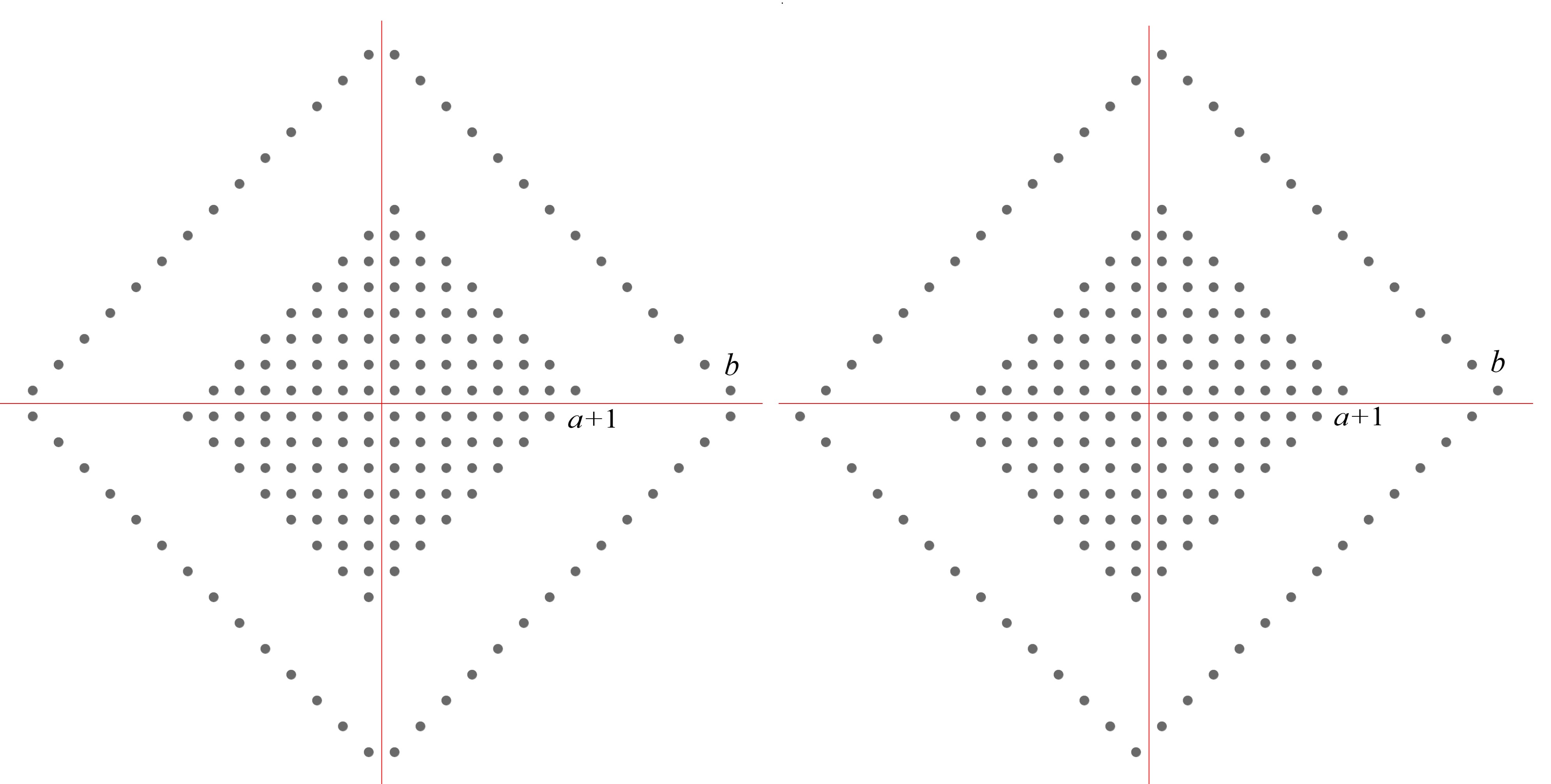} 
	\caption{$a+2\leq b\leq 2a$}
	\label{induct-2}
\end{figure}

Let $a\in\mathbb{N}$, $a\geq 4$. Consider four types of configurations given in Figures \ref{induct-1} and \ref{induct-2} (letters on the graph refer to $x$-coordinate; in all specific examples, $a=7$). All of them pass tests $\texttt{0}$, $\texttt{A}$, $\texttt{B}$, $\texttt{C}$. In particular, Figure \ref{induct-1} covers all even numbers in the range $[2a^2+6a+4,2a^2+10a]$. Figure \ref{induct-2} covers all even numbers in the range $[2a^2+8a+10,2a^2+12a+12]$.
Now, recall that Figure \ref{ger-p} (right) gives a particular example of inseparable set for $n=22$. This give the next result. 
\begin{prop}All even numbers with $c(n)=0$ are given by $n=2,6,10$.
	\end{prop}

\subsection{Odd numbers with $c(n)>0$} 
\begin{figure}
	\includegraphics[scale=0.22]{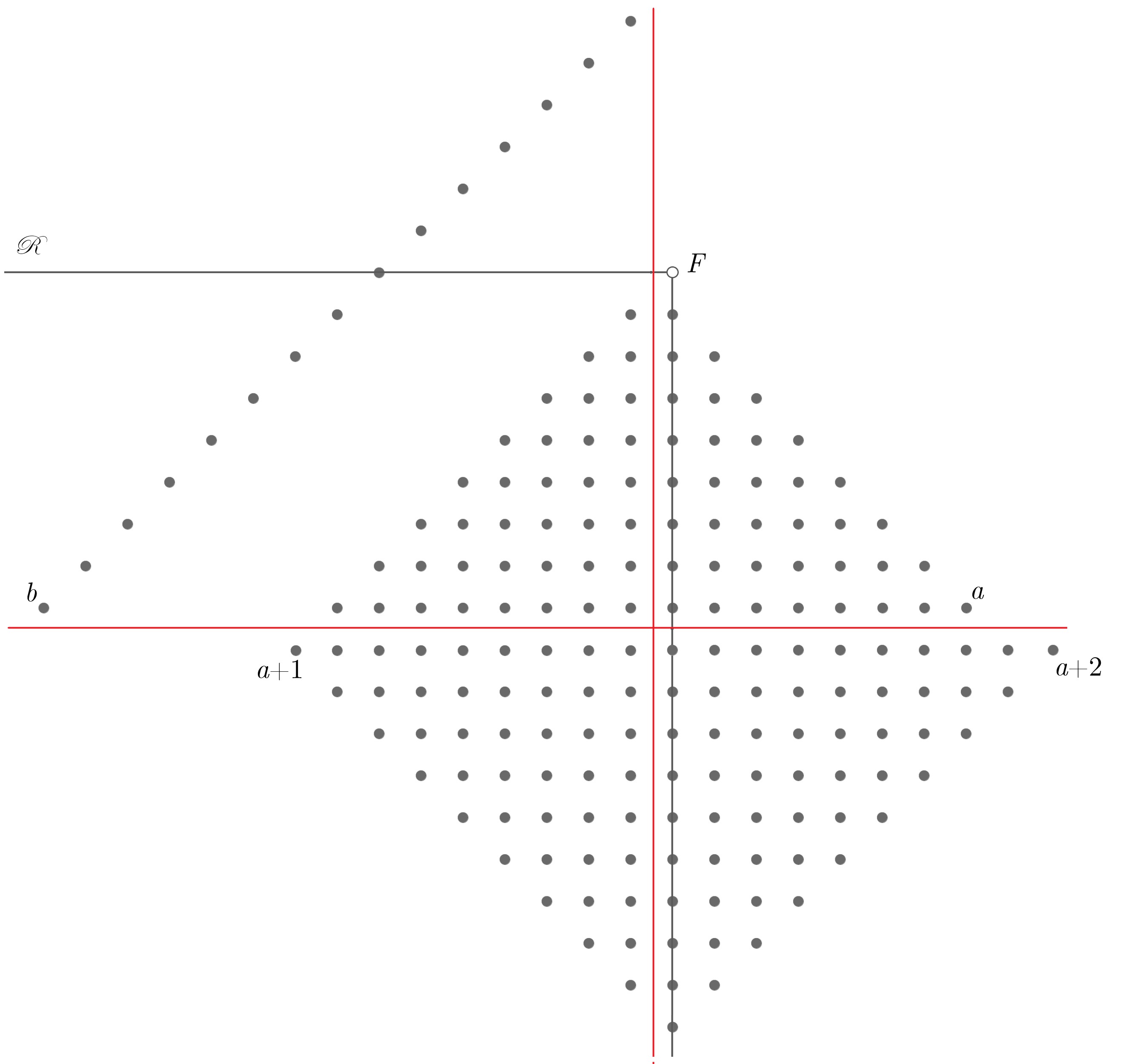} 
	\caption{The first family (the simplest one) of inseparable sets. Here $a\geq 2$; $b-a$ is odd; $a+3\leq b\leq 3a-1$ (type $[a,a+1,a+1,a+2]$, variation $2$).}
	\label{ger-ind}
\end{figure}

\begin{prop}Let $K(n)=\{i\in 2\mathbb{N}-1,i\leq n: c(i)=0\}$. Then
	\begin{eqnarray*}
		\lim\limits_{n\rightarrow\infty}\frac{K(n)}{n/2}=0.
		\end{eqnarray*}
	\label{visi-nel}
\end{prop}
In order to show this, we will present two carefully chosen families of inseparable $n$-sets. Let $a,b\in\mathbb{N}$, $b>a\geq 2$, $b-a$ being odd. Consider the construction in Figure \ref{ger-ind}. Here $n=2a^2+5a+b+4$. For a quartering to be proper, we must require $b\geq a+3$. Next, consider the path $\mathscr{R}$. Here $r(\mathscr{R})=a^2+3a+1$, and $\ell(\mathscr{R})=a^2+b$ points. If $b\leq 3a-1$, test $\mathtt{A}$ (subcase \textbf{3ii}) is passed. We easily calculate that test $\mathtt{A}$ then is passed in general. The construction accounts for all odd numbers in the interval $[2a^2+6a+7,2a^2+8a+3]$. For $a=2$ this reduces a singleton $n=27$, and for $a=3$ this give two numbers $n=43$ and $n=45$.\\

\begin{figure}
	\includegraphics[scale=0.26]{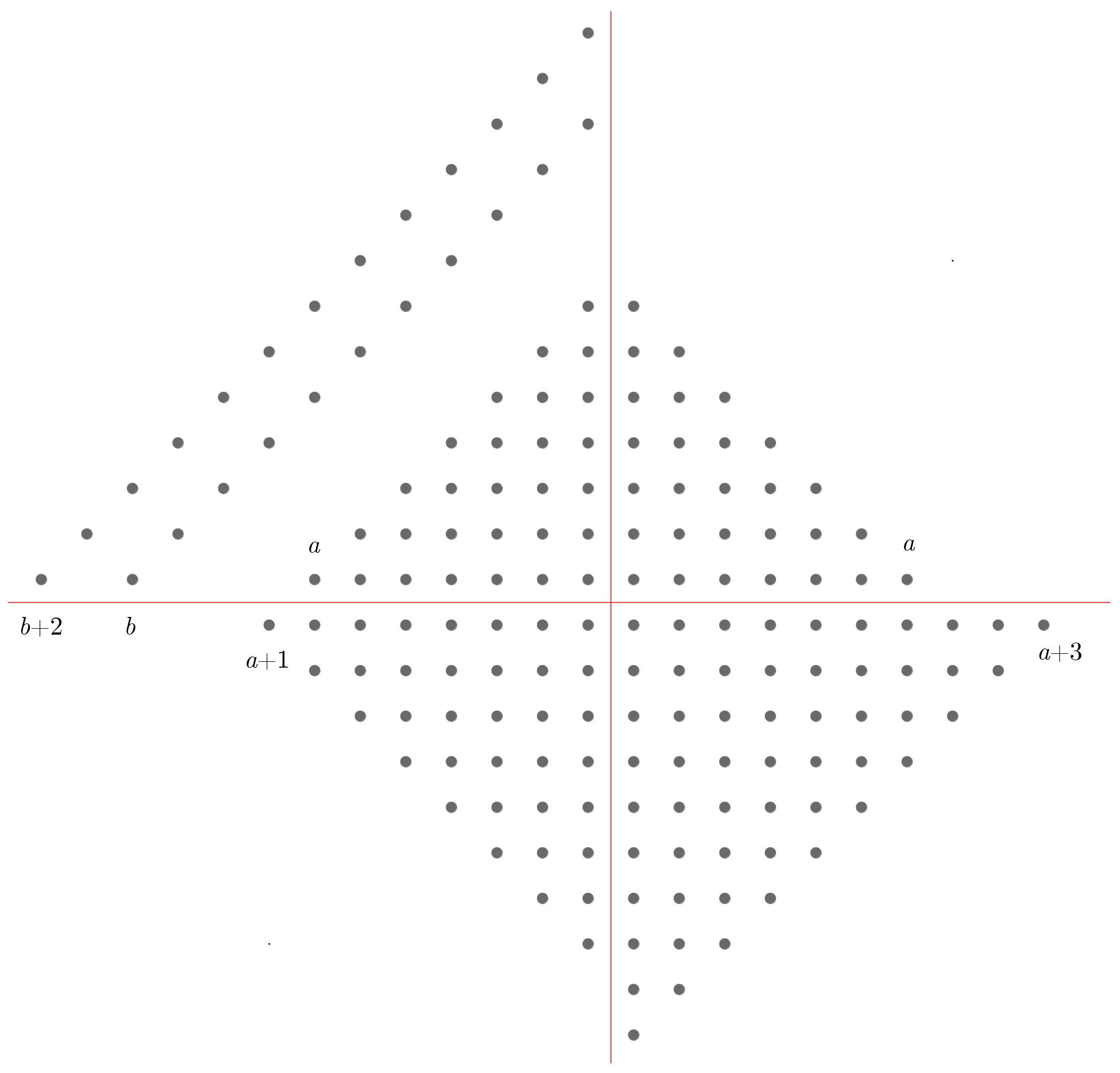} 
	\caption{The second family of inseparable sets. Here $a\geq 5$; $a+3\leq b\leq 2a-2$ (type $[a,a+2,a+1,a+3]$, variation $3$).}
	\label{ger-ind3}
\end{figure}

In the same manner, family in Figure \ref{ger-ind3} covers all odd integers in the range $[2a^2+8a+15, 2a^2+10a+5]$. Both these families cover all odd positive integers, save the ones contained in the intervals $[2a^2+8a+5,2a^2+8a+13]$ and $[2a^2+10a+7,2a^2+10a+13]$, $a\geq 5$. These (among odd ones) have a natural density $0$.

\subsection{Expectations}
\label{expect}
 What kind of sequences $c(n)$ and $\hat{c}(n)$ are expected to be? \\

Let us get back to Proposition \ref{prop-struc}. If we were to count only quadruples $(Y_{1},Y_{2},Y_{3},Y_{4})$ of positive integers which add up to $n$, the total number would be $\binom{n-1}{3}$ (sequence A000292 in \cite{oeis}). This is a number of \emph{compositions} (ordered partitions) of $n$ into $4$ positive parts. Let us now count the same compositions, only up to $\mathbb{D}_{4}$ symmetry (i.e. reversing order and cyclic shifts are allowed). Denote this number by $t(n)$. For example, $t(8)=8$, and these compositions are
\begin{eqnarray*}
	(1,1,1,5),\, (1,1,2,4),\, (1,1,3,3),\, (1,2,1,4),\,(1,3,1,3),\,(1,2,2,3),\,(1,2,3,2),\,(2,2,2,2).
\end{eqnarray*}
In general, if we count separately compositions having exactly $s$ summands equal to $1$ ($s=0,1,2$ or $3$), we arrive to the recurrence 
\begin{eqnarray*}
	t(n)=t(n-4)+\frac{1}{2}\binom{n-5}{2}+\frac{1}{2}\bigg{\lfloor}\frac{n-5}{2}\bigg{\rfloor}+2\bigg{\lfloor}\frac{n-4}{2}\bigg{\rfloor}+1,\quad n\geq 5.
\end{eqnarray*}
Thus, we get the sequence 
\begin{eqnarray*}
	0, 0, 0, 1, 1, 3, 4, 8, 10, 16, 20, 29, 35, 47, 56, 72, 84, 104, 120, 145,\ldots.
\end{eqnarray*}
This is A005232 in \cite{oeis} shifted by $4$: $t(n)=A005232(n-4)$. The generating function
\begin{eqnarray*}
	\sum\limits_{n=1}^{\infty}t(n)z^{n}=\frac{z^4(1-z+z^2)}{(1-z)^2(1-z^2)(1-z^4)}.
\end{eqnarray*}
But now, let us take a look at the Table 2. These are the two examples of inseparable sets for $n=44$:
\begin{eqnarray*}
[(1,2, \mathbf{3}, \mathbf{6})(1, 2, 3, 4)(1, 2, 3, 6)(1, 2, 3, 4)],\,[(1, 2, \mathbf{4}, \mathbf{5})(1, 2, 3, 4)(1, 2, 3, 6)(1, 2, 3, 4)].
\end{eqnarray*}
They are of the forms $[L_{1}L_{2}L_{3}L_{4}]$ and $[L_{1}'L_{2}L_{3}L_{4}]$. The only difference are numbers marked in bold. In particular, this is possible due to equality $3+6=4+5$. Though being a trivial observation, it shows that the function $q(n)$ of \emph{partitions into distinct parts} (A000009 in \cite{oeis}) manifests in the structure of both sequences $c(n)$ and $\hat{c}(n)$. All these question will be treated in detail in \cite{antra-dalis}.

\section{Supplements}
\begin{center}
	\setlength{\shadowsize}{2pt}\shadowbox{\texttt{Test B}}
\end{center}
\begin{verbatim}
B:=proc(j,i::integer)
local Tt, t, wij, g, h, lik, v, Ls, Rs, yra;
global L,N,n,m;
yra:=true;
t:=2;
if L[i][1]>1 then
if abs(Y[i]-Y[modp(i+1,4)+1])>
Y[modp(i,4)+1]+Y[modp(i-2,4)+1]-N[modp(i,4)+1]-N[modp(i-2,4)+1] 
then Ls:=1: Rs:=0: Tt:=0: else Ls:=0: Rs:=1: Tt:=0: end if: end if:
if L[i][1]=1 then 
while member(t,L[i])=member(t-1,L[j]) and t<=max(m[i],m[j]) do t:=t+1: end do:
if member(t,L[i]) then Tt:=N[modp(j+1,4)+1]+N[j]+2: end if:
if member(t-1,L[j])then Tt:=N[modp(j+1,4)+1]+N[j]: end if:
if t=max(m[i],m[j])+1 then yra:=false: end if:
if yra then
wij:=t: 
if L[i][N[i]]<=wij then h:=N[i] else
h:=1:
while (L[i][h]<=wij) do h:=h+1: end do:
h:=h-1:
end if:
g:=Tt-h-N[modp(j+1,4)+1]:
lik:=0:
for v from 1 to N[j] do if L[j][v]>wij then lik:=lik+(L[j][v]-wij): end if: end do:
Ls:=Y[modp(i+1,4)+1]+Y[j]-g-lik: 
Rs:=n-Tt-Ls:
else Ls:=1: Rs:=0: end if: end if:
[Ls,Rs,yra,Tt]:
end proc:

\end{verbatim}

\begin{center}
	\setlength{\shadowsize}{2pt}\shadowbox{\texttt{Test C}}
\end{center}
\begin{verbatim}
	C:=proc(i::integer)
	local s, t,lik, Tt, v, Rs, Ls, Ger;
	global L,N,n,m;
	s:=1;
	while member(s,L[i]) do s:=s+1: end do:
	Ger:=true:
	for t from 1 to s-2 do lik:=0: Tt:=0:
	if member(s-t,L[modp(i,4)+1]) then Tt:=Tt+1: end if:
	for v from 1 to N[modp(i,4)+1] do if L[modp(i,4)+1][v]>s-t then
	lik:=lik+(L[modp(i,4)+1][v]-s+t): Tt:=Tt+1: end if: end do:
	if member(t+1,L[modp(i-2,4)+1]) then Tt:=Tt+1: end if:
	for v from 1 to N[modp(i-2,4)+1] do if L[modp(i-2,4)+1][v]>t+1 then
	lik:=lik+(L[modp(i-2,4)+1][v]-t-1): Tt:=Tt+1: end if: end do:
	Rs:=lik+Y[i]-(t+1)*(s-t)+1: Tt:=Tt+s-1: Ls:=n-Tt-Rs:
	Ger:=Ger and (Ls>Rs):
	end do:
	Ger:
	end proc:	
\end{verbatim}

\begin{figure}
		\includegraphics[scale=0.78]{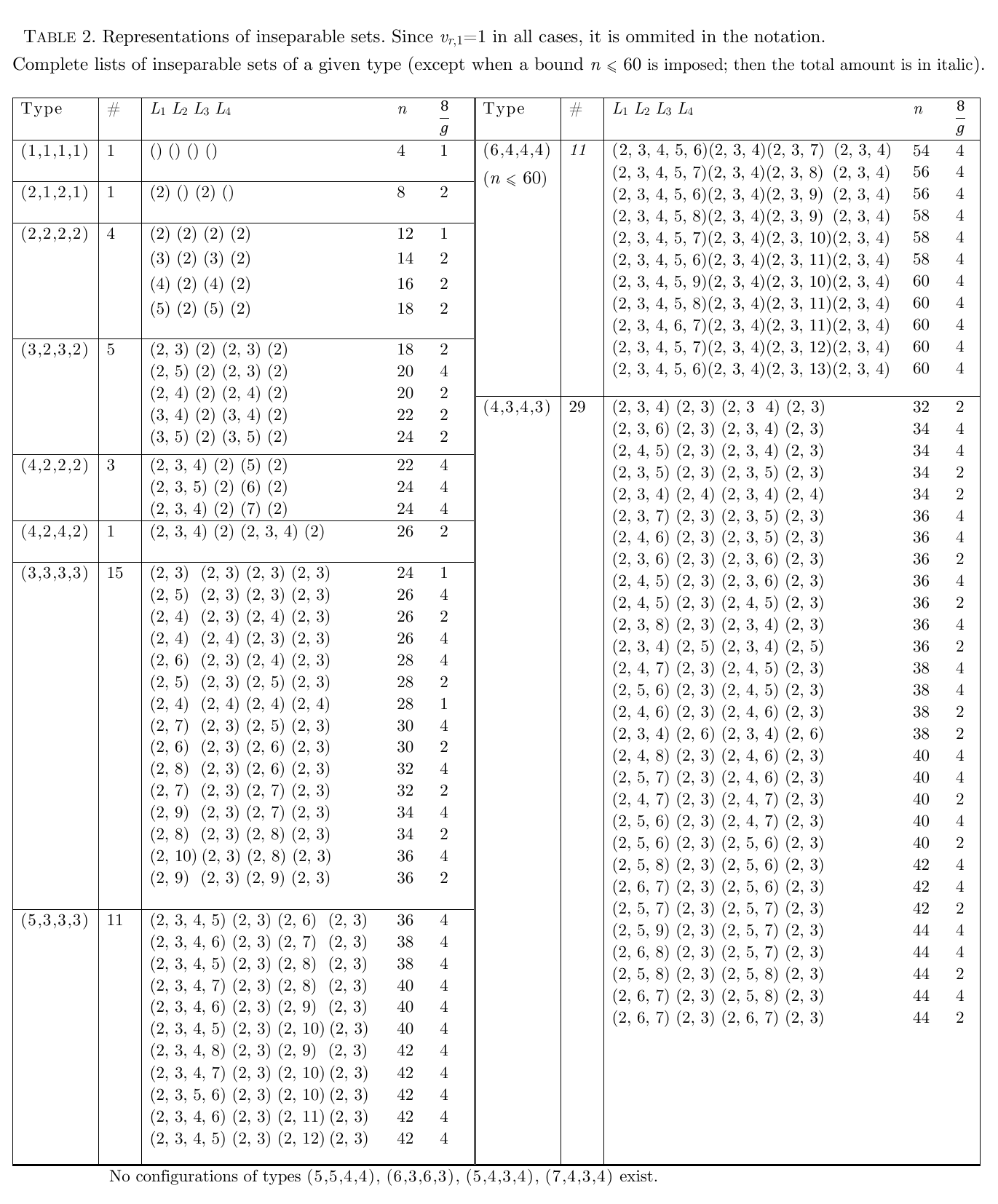}
\end{figure}  
\begin{figure}
		\includegraphics[scale=0.78]{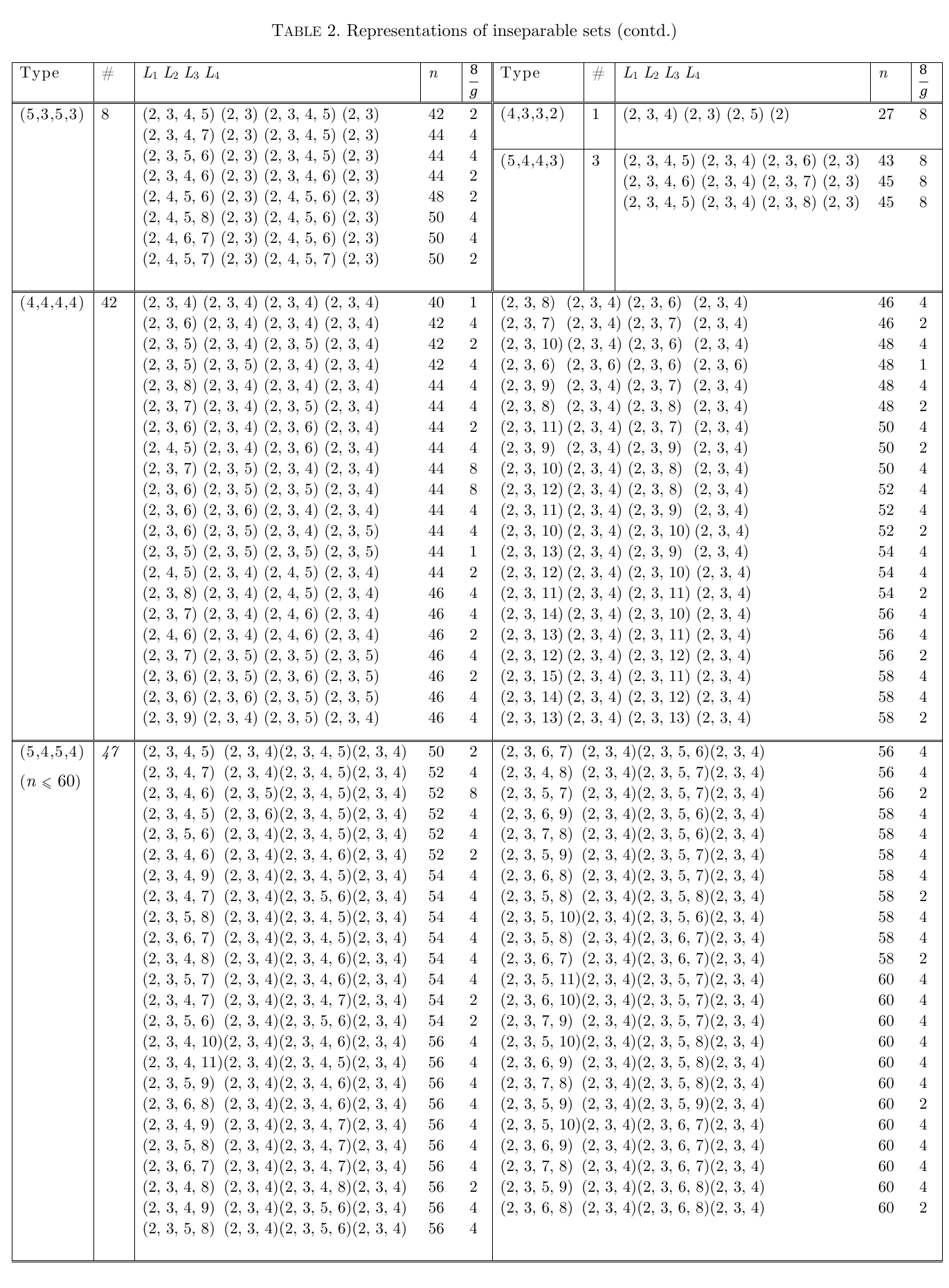}
\end{figure}

\end{document}